\documentclass{amsart}
\usepackage{amsfonts}

\newtheorem{theorem}{Theorem}
\theoremstyle{plain}

\newtheorem{corollary}{Corollary}

\newtheorem{proposition}{Proposition}
\newtheorem{remark}{Remark}

\numberwithin{equation}{section}
\input{tcilatex}

\begin{document}
\title[Reverses of the Continuous Triangle Inequality]{Reverses of the
Continuous Triangle Inequality for Bochner Integral of Vector-Valued
Functions in Hilbert Spaces}
\author{Sever S. Dragomir}
\address{School of Computer Science and Mathematics\\
Victoria University of Technology\\
PO Box 14428, MCMC 8001\\
Victoria, Australia.}
\email{sever@csm.vu.edu.au}
\urladdr{http://rgmia.vu.edu.au/SSDragomirWeb.html}
\date{April 16, 2004.}

\begin{abstract}
Some reverses of the continuous triangle inequality for Bochner integral of
vector-valued functions in Hilbert spaces are given. Applications for
complex-valued functions are provided as well.
\end{abstract}

\maketitle

\section{Introduction}

Let $f:\left[ a,b\right] \rightarrow \mathbb{K}$, $\mathbb{K}=\mathbb{C}$ or 
$\mathbb{R}$ be a Lebesgue integrable function. The following inequality,
which is the continuous version of the \textit{triangle inequality}%
\begin{equation}
\left\vert \int_{a}^{b}f\left( x\right) dx\right\vert \leq
\int_{a}^{b}\left\vert f\left( x\right) \right\vert dx,  \label{1.1}
\end{equation}%
plays a fundamental role in Mathematical Analysis and its applications.

It appears, see \cite[p. 492]{MPF}, that the first reverse inequality for (%
\ref{1.1}) was obtained by J. Karamata in his book from 1949, \cite{K}. It
can be stated as%
\begin{equation}
\cos \theta \int_{a}^{b}\left\vert f\left( x\right) \right\vert dx\leq
\left\vert \int_{a}^{b}f\left( x\right) dx\right\vert  \label{1.2}
\end{equation}%
provided%
\begin{equation*}
-\theta \leq \arg f\left( x\right) \leq \theta ,\ \ x\in \left[ a,b\right]
\end{equation*}%
for given $\theta \in \left( 0,\frac{\pi }{2}\right) .$

This integral inequality is the continuous version of a reverse inequality
for the generalised triangle inequality%
\begin{equation}
\cos \theta \sum_{i=1}^{n}\left\vert z_{i}\right\vert \leq \left\vert
\sum_{i=1}^{n}z_{i}\right\vert ,  \label{1.3}
\end{equation}%
provided%
\begin{equation*}
a-\theta \leq \arg \left( z_{i}\right) \leq a+\theta ,\ \ \text{for \ }i\in
\left\{ 1,\dots ,n\right\} ,
\end{equation*}%
where $a\in \mathbb{R}$ and $\theta \in \left( 0,\frac{\pi }{2}\right) ,$
which, as pointed out in \cite[p. 492]{MPF}, was first discovered by M.
Petrovich in 1917, \cite{P}, and, subsequently rediscovered by other
authors, including J. Karamata \cite[p. 300 -- 301]{K}, H.S. Wilf \cite{W},
and in an equivalent form, by M. Marden \cite{M}.

The first to consider the problem in the more general case of Hilbert and
Banach spaces, were J.B. Diaz and F.T. Metcalf \cite{DM} who showed that, in
an inner product space $H$ over the real or complex number field, the
following reverse of the triangle inequality holds%
\begin{equation}
r\sum_{i=1}^{n}\left\Vert x_{i}\right\Vert \leq \left\Vert
\sum_{i=1}^{n}x_{i}\right\Vert ,  \label{1.4}
\end{equation}%
provided%
\begin{equation*}
0\leq r\leq \frac{\func{Re}\left\langle x_{i},a\right\rangle }{\left\Vert
x_{i}\right\Vert },\ \ \ \ \ i\in \left\{ 1,\dots ,n\right\} ,
\end{equation*}%
and $a\in H$ is a unit vector, i.e., $\left\Vert a\right\Vert =1.$ The case
of equality holds in (\ref{1.4}) if and only if%
\begin{equation}
\sum_{i=1}^{n}x_{i}=r\left( \sum_{i=1}^{n}\left\Vert x_{i}\right\Vert
\right) a.  \label{1.5}
\end{equation}

The main aim of this paper is to point out some reverses of the triangle
inequality for Bochner integrable functions $f$ with values in Hilbert
spaces and defined on a compact interval $\left[ a,b\right] \subset \mathbb{R%
}$. Applications for Lebesgue integrable complex-valued functions are
provided as well.

\section{Reverses for a Unit Vector}

We recall that $f\in L\left( \left[ a,b\right] ;H\right) ,$ the space of
Bochner integrable functions with values in a Hilbert space $H,$ if and only
if $f:\left[ a,b\right] \rightarrow H$ is Bochner measurable on $\left[ a,b%
\right] $ and the Lebesgue integral $\int_{a}^{b}\left\Vert f\left( t\right)
\right\Vert dt$ is finite.

The following result holds:

\begin{theorem}
\label{t2.1}If $f\in L\left( \left[ a,b\right] ;H\right) $ is such that
there exists a constant $K\geq 1$ and a vector $e\in H,$ $\left\Vert
e\right\Vert =1$ with%
\begin{equation}
\left\Vert f\left( t\right) \right\Vert \leq K\func{Re}\left\langle f\left(
t\right) ,e\right\rangle \ \ \ \text{for a.e. }t\in \left[ a,b\right] ,
\label{2.1}
\end{equation}%
then we have the inequality:%
\begin{equation}
\int_{a}^{b}\left\Vert f\left( t\right) \right\Vert dt\leq K\left\Vert
\int_{a}^{b}f\left( t\right) dt\right\Vert .  \label{2.2}
\end{equation}%
The case of equality holds in (\ref{2.2}) if and only if%
\begin{equation}
\int_{a}^{b}f\left( t\right) dt=\frac{1}{K}\left( \int_{a}^{b}\left\Vert
f\left( t\right) \right\Vert dt\right) e.  \label{2.3}
\end{equation}
\end{theorem}

\begin{proof}
By the Schwarz inequality in inner product spaces, we have%
\begin{align}
\left\Vert \int_{a}^{b}f\left( t\right) dt\right\Vert & =\left\Vert
\int_{a}^{b}f\left( t\right) dt\right\Vert \left\Vert e\right\Vert
\label{2.4} \\
& \geq \left\vert \left\langle \int_{a}^{b}f\left( t\right)
dt,e\right\rangle \right\vert \geq \left\vert \func{Re}\left\langle
\int_{a}^{b}f\left( t\right) dt,e\right\rangle \right\vert  \notag \\
& \geq \func{Re}\left\langle \int_{a}^{b}f\left( t\right) dt,e\right\rangle
=\int_{a}^{b}\func{Re}\left\langle f\left( t\right) ,e\right\rangle dt. 
\notag
\end{align}%
From the condition (\ref{2.1}), on integrating over $\left[ a,b\right] ,$ we
deduce%
\begin{equation}
\int_{a}^{b}\func{Re}\left\langle f\left( t\right) ,e\right\rangle dt\geq 
\frac{1}{K}\int_{a}^{b}\left\Vert f\left( t\right) \right\Vert dt,
\label{2.5}
\end{equation}%
and thus, on making use of (\ref{2.4}) and (\ref{2.5}), we obtain the
desired inequality (\ref{2.2}).

If (\ref{2.3}) holds true, then, obviously%
\begin{equation*}
K\left\Vert \int_{a}^{b}f\left( t\right) dt\right\Vert =\left\Vert
e\right\Vert \int_{a}^{b}\left\Vert f\left( t\right) \right\Vert
dt=\int_{a}^{b}\left\Vert f\left( t\right) \right\Vert dt,
\end{equation*}%
showing that (\ref{2.2}) holds with equality.

If we assume that the equality holds in (\ref{2.2}), then by the argument
provided at the beginning of our proof, we must have equality in each of the
inequalities from (\ref{2.4}) and (\ref{2.5}).

Observe that in Schwarz's inequality $\left\Vert x\right\Vert \left\Vert
y\right\Vert \geq \func{Re}\left\langle x,y\right\rangle ,$ $x,y\in H,$ the
case of equality holds if and only if there exists a positive scalar $\mu $
such that $x=\mu e.$ Therefore, equality holds in the first inequality in (%
\ref{2.4}) iff $\int_{a}^{b}f\left( t\right) dt=\lambda e$, with $\lambda
\geq 0$ $.$

If we assume that a strict inequality holds in (\ref{2.1}) for a.e. $t\in %
\left[ a,b\right] ,$ then $\int_{a}^{b}\left\Vert f\left( t\right)
\right\Vert dt<K\int_{a}^{b}\func{Re}\left\langle f\left( t\right)
,e\right\rangle dt,$ and by (\ref{2.4}) we deduce a strict inequality in (%
\ref{2.2}), which contradicts the assumption. Thus, we must have $\left\Vert
f\left( t\right) \right\Vert =K\func{Re}\left\langle f\left( t\right)
,e\right\rangle $ for a.e. $t\in \left[ a,b\right] .$

If we integrate this equality, we deduce%
\begin{align*}
\int_{a}^{b}\left\Vert f\left( t\right) \right\Vert dt& =K\int_{a}^{b}\func{%
Re}\left\langle f\left( t\right) ,e\right\rangle dt=K\func{Re}\left\langle
\int_{a}^{b}f\left( t\right) dt,e\right\rangle \\
& =K\func{Re}\left\langle \lambda e,e\right\rangle =\lambda K
\end{align*}%
giving%
\begin{equation*}
\lambda =\frac{1}{K}\int_{a}^{b}\left\Vert f\left( t\right) \right\Vert dt,
\end{equation*}%
and thus the equality (\ref{2.3}) is necessary.

This completes the proof.
\end{proof}

A more appropriate result from an applications point of view is perhaps the
following result.

\begin{corollary}
\label{c2.2}Let $e$ be a unit vector in the Hilbert space $\left(
H;\left\langle \cdot ,\cdot \right\rangle \right) ,$ $\rho \in \left(
0,1\right) $ and $f\in L\left( \left[ a,b\right] ;H\right) $ so that%
\begin{equation}
\left\Vert f\left( t\right) -e\right\Vert \leq \rho \ \ \ \text{for a.e. }%
t\in \left[ a,b\right] .  \label{2.6}
\end{equation}%
Then we have the inequality%
\begin{equation}
\sqrt{1-\rho ^{2}}\int_{a}^{b}\left\Vert f\left( t\right) \right\Vert dt\leq
\left\Vert \int_{a}^{b}f\left( t\right) dt\right\Vert ,  \label{2.7}
\end{equation}%
with equality if and only if%
\begin{equation}
\int_{a}^{b}f\left( t\right) dt=\sqrt{1-\rho ^{2}}\left(
\int_{a}^{b}\left\Vert f\left( t\right) \right\Vert dt\right) \cdot e.
\label{2.8}
\end{equation}
\end{corollary}

\begin{proof}
From (\ref{2.6}), we have%
\begin{equation*}
\left\Vert f\left( t\right) \right\Vert ^{2}-2\func{Re}\left\langle f\left(
t\right) ,e\right\rangle +1\leq \rho ^{2},
\end{equation*}%
giving%
\begin{equation*}
\left\Vert f\left( t\right) \right\Vert ^{2}+1-\rho ^{2}\leq 2\func{Re}%
\left\langle f\left( t\right) ,e\right\rangle
\end{equation*}%
for a.e. $t\in \left[ a,b\right] .$

Dividing by $\sqrt{1-\rho ^{2}}>0,$ we deduce%
\begin{equation}
\frac{\left\Vert f\left( t\right) \right\Vert ^{2}}{\sqrt{1-\rho ^{2}}}+%
\sqrt{1-\rho ^{2}}\leq \frac{2\func{Re}\left\langle f\left( t\right)
,e\right\rangle }{\sqrt{1-\rho ^{2}}}  \label{2.9}
\end{equation}%
for a.e. $t\in \left[ a,b\right] .$

On the other hand, by the elementary inequality%
\begin{equation*}
\frac{p}{\alpha }+q\alpha \geq 2\sqrt{pq},\ \ \ p,q\geq 0,\ \alpha >0
\end{equation*}%
we have%
\begin{equation}
2\left\Vert f\left( t\right) \right\Vert \leq \frac{\left\Vert f\left(
t\right) \right\Vert ^{2}}{\sqrt{1-\rho ^{2}}}+\sqrt{1-\rho ^{2}}
\label{2.10}
\end{equation}%
for each $t\in \left[ a,b\right] .$

Making use of (\ref{2.9}) and (\ref{2.10}), we deduce%
\begin{equation*}
\left\Vert f\left( t\right) \right\Vert \leq \frac{1}{\sqrt{1-\rho ^{2}}}%
\func{Re}\left\langle f\left( t\right) ,e\right\rangle
\end{equation*}%
for a.e. $t\in \left[ a,b\right] .$

Applying Theorem \ref{t2.1} for $K=\frac{1}{\sqrt{1-\rho ^{2}}},$ we deduce
the desired inequality (\ref{2.7}).
\end{proof}

In the same spirit, we also have the following corollary.

\begin{corollary}
\label{c2.3}Let $e$ be a unit vector in $H$ and $M\geq m>0.$ If $f\in
L\left( \left[ a,b\right] ;H\right) $ is such that%
\begin{equation}
\func{Re}\left\langle Me-f\left( t\right) ,f\left( t\right) -me\right\rangle
\geq 0\ \ \ \text{for a.e. }t\in \left[ a,b\right] ,  \label{2.9a}
\end{equation}%
or, equivalently,%
\begin{equation}
\left\Vert f\left( t\right) -\frac{M+m}{2}e\right\Vert \leq \frac{1}{2}%
\left( M-m\right) \ \ \ \text{for a.e. }t\in \left[ a,b\right] ,
\label{2.10a}
\end{equation}%
then we have the inequality%
\begin{equation}
\frac{2\sqrt{mM}}{M+m}\int_{a}^{b}\left\Vert f\left( t\right) \right\Vert
dt\leq \left\Vert \int_{a}^{b}f\left( t\right) dt\right\Vert ,  \label{2.11}
\end{equation}%
or, equivalently,%
\begin{equation}
\left( 0\leq \right) \int_{a}^{b}\left\Vert f\left( t\right) \right\Vert
dt-\left\Vert \int_{a}^{b}f\left( t\right) dt\right\Vert \leq \frac{\left( 
\sqrt{M}-\sqrt{m}\right) ^{2}}{M+m}\left\Vert \int_{a}^{b}f\left( t\right)
dt\right\Vert .  \label{2.12}
\end{equation}%
The equality holds in (\ref{2.11}) (or in the second part of (\ref{2.12}))
if and only if%
\begin{equation}
\int_{a}^{b}f\left( t\right) dt=\frac{2\sqrt{mM}}{M+m}\left(
\int_{a}^{b}\left\Vert f\left( t\right) \right\Vert dt\right) e.
\label{2.13}
\end{equation}
\end{corollary}

\begin{proof}
Firstly, we remark that if $x,z,Z\in H,$ then the following statements are
equivalent

\begin{enumerate}
\item[(i)] $\func{Re}\left\langle Z-x,x-z\right\rangle \geq 0$

and

\item[(ii)] $\left\Vert x-\frac{Z+z}{2}\right\Vert \leq \frac{1}{2}%
\left\Vert Z-z\right\Vert .$
\end{enumerate}

Using this fact, we may simply realise that (\ref{2.9}) and (\ref{2.10}) are
equivalent.

Now, from (\ref{2.9}), we obtain%
\begin{equation*}
\left\Vert f\left( t\right) \right\Vert ^{2}+mM\leq \left( M+m\right) \func{%
Re}\left\langle f\left( t\right) ,e\right\rangle
\end{equation*}%
for a.e. $t\in \left[ a,b\right] .$ Dividing this inequality with $\sqrt{mM}%
>0,$ we deduce the following inequality that will be used in the sequel%
\begin{equation}
\frac{\left\Vert f\left( t\right) \right\Vert ^{2}}{\sqrt{mM}}+\sqrt{mM}\leq 
\frac{M+m}{\sqrt{mM}}\func{Re}\left\langle f\left( t\right) ,e\right\rangle
\label{2.14}
\end{equation}%
for a.e. $t\in \left[ a,b\right] .$

On the other hand%
\begin{equation}
2\left\Vert f\left( t\right) \right\Vert \leq \frac{\left\Vert f\left(
t\right) \right\Vert ^{2}}{\sqrt{mM}}+\sqrt{mM},  \label{2.15}
\end{equation}%
for any $t\in \left[ a,b\right] .$

Utilising (\ref{2.14}) and (\ref{2.15}), we may conclude with the following
inequality 
\begin{equation*}
\left\Vert f\left( t\right) \right\Vert \leq \frac{M+m}{2\sqrt{mM}}\func{Re}%
\left\langle f\left( t\right) ,e\right\rangle ,
\end{equation*}%
for a.e. $t\in \left[ a,b\right] .$

Applying Theorem \ref{t2.1} for the constant $K:=\frac{m+M}{2\sqrt{mM}}\geq
1,$ we deduce the desired result.
\end{proof}

\section{Reverses for Orthornormal Families of Vectors}

The following result for orthornormal vectors in $H$ holds.

\begin{theorem}
\label{t3.1}Let $\left\{ e_{1},\dots ,e_{n}\right\} $ be a family of
orthornormal vectors in $H,$ $k_{i}\geq 0,$ $i\in \left\{ 1,\dots ,n\right\} 
$ and $f\in L\left( \left[ a,b\right] ;H\right) $ such that%
\begin{equation}
k_{i}\left\Vert f\left( t\right) \right\Vert \leq \func{Re}\left\langle
f\left( t\right) ,e_{i}\right\rangle  \label{3.1}
\end{equation}%
for each $i\in \left\{ 1,\dots ,n\right\} $ and for a.e. $t\in \left[ a,b%
\right] .$

Then%
\begin{equation}
\left( \sum_{i=1}^{n}k_{i}^{2}\right) ^{\frac{1}{2}}\int_{a}^{b}\left\Vert
f\left( t\right) \right\Vert dt\leq \left\Vert \int_{a}^{b}f\left( t\right)
dt\right\Vert ,  \label{3.2}
\end{equation}%
where equality holds if and only if%
\begin{equation}
\int_{a}^{b}f\left( t\right) dt=\left( \int_{a}^{b}\left\Vert f\left(
t\right) \right\Vert dt\right) \sum_{i=1}^{n}k_{i}e_{i}.  \label{3.3}
\end{equation}
\end{theorem}

\begin{proof}
By Bessel's inequality applied for $\int_{a}^{b}f\left( t\right) dt$ and the
orthornormal vectors $\left\{ e_{1},\dots ,e_{n}\right\} ,$ we have%
\begin{align}
\left\Vert \int_{a}^{b}f\left( t\right) dt\right\Vert ^{2}& \geq
\sum_{i=1}^{n}\left\vert \left\langle \int_{a}^{b}f\left( t\right)
dt,e_{i}\right\rangle \right\vert ^{2}  \label{3.4} \\
& \geq \sum_{i=1}^{n}\left[ \func{Re}\left\langle \int_{a}^{b}f\left(
t\right) dt,e_{i}\right\rangle \right] ^{2}  \notag \\
& =\sum_{i=1}^{n}\left[ \int_{a}^{b}\func{Re}\left\langle f\left( t\right)
,e_{i}\right\rangle dt\right] ^{2}.  \notag
\end{align}%
Integrating (\ref{3.1}), we get for each $i\in \left\{ 1,\dots ,n\right\} $%
\begin{equation*}
0\leq k_{i}\int_{a}^{b}\left\Vert f\left( t\right) \right\Vert dt\leq
\int_{a}^{b}\func{Re}\left\langle f\left( t\right) ,e_{i}\right\rangle dt,
\end{equation*}%
implying%
\begin{equation}
\sum_{i=1}^{n}\left[ \int_{a}^{b}\func{Re}\left\langle f\left( t\right)
,e_{i}\right\rangle dt\right] ^{2}\geq \sum_{i=1}^{n}k_{i}^{2}\left(
\int_{a}^{b}\left\Vert f\left( t\right) \right\Vert dt\right) ^{2}.
\label{3.5}
\end{equation}%
On making use of (\ref{3.4}) and (\ref{3.5}), we deduce%
\begin{equation*}
\left\Vert \int_{a}^{b}f\left( t\right) dt\right\Vert ^{2}\geq
\sum_{i=1}^{n}k_{i}^{2}\left( \int_{a}^{b}\left\Vert f\left( t\right)
\right\Vert dt\right) ^{2},
\end{equation*}%
which is clearly equivalent to (\ref{3.2}).

If (\ref{3.3}) holds true, then%
\begin{align*}
\left\Vert \int_{a}^{b}f\left( t\right) dt\right\Vert & =\left(
\int_{a}^{b}\left\Vert f\left( t\right) \right\Vert dt\right) \left\Vert
\sum_{i=1}^{n}k_{i}e_{i}\right\Vert \\
& =\left( \int_{a}^{b}\left\Vert f\left( t\right) \right\Vert dt\right) 
\left[ \sum_{i=1}^{n}k_{i}^{2}\left\Vert e_{i}\right\Vert ^{2}\right] ^{%
\frac{1}{2}} \\
& =\left( \sum_{i=1}^{n}k_{i}^{2}\right) ^{\frac{1}{2}}\int_{a}^{b}\left%
\Vert f\left( t\right) \right\Vert dt,
\end{align*}%
showing that (\ref{3.2}) holds with equality.

Now, suppose that there is an $i_{0}\in \left\{ 1,\dots ,n\right\} $ for
which%
\begin{equation*}
k_{i_{0}}\left\Vert f\left( t\right) \right\Vert <\func{Re}\left\langle
f\left( t\right) ,e_{i_{0}}\right\rangle
\end{equation*}%
for a.e. $t\in \left[ a,b\right] .$ Then obviously%
\begin{equation*}
k_{i_{0}}\int_{a}^{b}\left\Vert f\left( t\right) \right\Vert dt<\int_{a}^{b}%
\func{Re}\left\langle f\left( t\right) ,e_{i_{0}}\right\rangle dt,
\end{equation*}%
and using the argument given above, we deduce%
\begin{equation*}
\left( \sum_{i=1}^{n}k_{i}^{2}\right) ^{\frac{1}{2}}\int_{a}^{b}\left\Vert
f\left( t\right) \right\Vert dt<\left\Vert \int_{a}^{b}f\left( t\right)
dt\right\Vert .
\end{equation*}%
Therefore, if the equality holds in (\ref{3.2}), we must have%
\begin{equation}
k_{i}\left\Vert f\left( t\right) \right\Vert =\func{Re}\left\langle f\left(
t\right) ,e_{i}\right\rangle  \label{3.6}
\end{equation}%
for each $i\in \left\{ 1,\dots ,n\right\} $ and a.e. $t\in \left[ a,b\right]
.$

Also, if the equality holds in (\ref{3.2}), then we must have equality in
all inequalities (\ref{3.4}), this means that%
\begin{equation}
\int_{a}^{b}f\left( t\right) dt=\sum_{i=1}^{n}\left\langle
\int_{a}^{b}f\left( t\right) dt,e_{i}\right\rangle e_{i}  \label{3.7}
\end{equation}%
and 
\begin{equation}
\func{Im}\left\langle \int_{a}^{b}f\left( t\right) dt,e_{i}\right\rangle =0%
\text{ \ for each \ }i\in \left\{ 1,\dots ,n\right\} .  \label{3.8}
\end{equation}%
Using (\ref{3.6}) and (\ref{3.8}) in (\ref{3.7}), we deduce%
\begin{align*}
\int_{a}^{b}f\left( t\right) dt& =\sum_{i=1}^{n}\func{Re}\left\langle
\int_{a}^{b}f\left( t\right) dt,e_{i}\right\rangle e_{i} \\
& =\sum_{i=1}^{n}\int_{a}^{b}\func{Re}\left\langle f\left( t\right)
,e_{i}\right\rangle e_{i}dt \\
& =\sum_{i=1}^{n}\left( \int_{a}^{b}\left\Vert f\left( t\right) \right\Vert
dt\right) k_{i}e_{i} \\
& =\int_{a}^{b}\left\Vert f\left( t\right) \right\Vert
dt\sum_{i=1}^{n}k_{i}e_{i},
\end{align*}%
and the condition (\ref{3.3}) is necessary.

This completes the proof.
\end{proof}

The following two corollaries are of interest.

\begin{corollary}
\label{c3.2}Let $\left\{ e_{1},\dots ,e_{n}\right\} $ be a family of
orthornormal vectors in $H,$ $\rho _{i}\in \left( 0,1\right) ,$ $i\in
\left\{ 1,\dots ,n\right\} $ and $f\in L\left( \left[ a,b\right] ;H\right) $
such that:%
\begin{equation}
\left\Vert f\left( t\right) -e_{i}\right\Vert \leq \rho _{i}\text{ \ for \ }%
i\in \left\{ 1,\dots ,n\right\} \text{ \ and a.e. }t\in \left[ a,b\right] .
\label{3.9}
\end{equation}%
Then we have the inequality%
\begin{equation*}
\left( n-\sum_{i=1}^{n}\rho _{i}^{2}\right) ^{\frac{1}{2}}\int_{a}^{b}\left%
\Vert f\left( t\right) \right\Vert dt\leq \left\Vert \int_{a}^{b}f\left(
t\right) dt\right\Vert ,
\end{equation*}%
with equality if and only if%
\begin{equation*}
\int_{a}^{b}f\left( t\right) dt=\int_{a}^{b}\left\Vert f\left( t\right)
\right\Vert dt\left( \sum_{i=1}^{n}\left( 1-\rho _{i}^{2}\right)
^{1/2}e_{i}\right) .
\end{equation*}
\end{corollary}

\begin{proof}
From the proof of Theorem \ref{t2.1}, we know that (\ref{3.3}) implies the
inequality%
\begin{equation*}
\sqrt{1-\rho _{i}^{2}}\left\Vert f\left( t\right) \right\Vert \leq \func{Re}%
\left\langle f\left( t\right) ,e_{i}\right\rangle ,\ \ \ i\in \left\{
1,\dots ,n\right\} ,\text{ \ for a.e. }t\in \left[ a,b\right] .
\end{equation*}%
Now, applying Theorem \ref{t3.1} for $k_{i}:=\sqrt{1-\rho _{i}^{2}},$ $i\in
\left\{ 1,\dots ,n\right\} $, we deduce the desired result.
\end{proof}

\begin{corollary}
\label{c3.3}Let $\left\{ e_{1},\dots ,e_{n}\right\} $ be a family of
orthornormal vectors in $H,$ $M_{i}\geq m_{i}>0,$ $i\in \left\{ 1,\dots
,n\right\} $ and $f\in L\left( \left[ a,b\right] ;H\right) $ such that%
\begin{equation}
\func{Re}\left\langle M_{i}e_{i}-f\left( t\right) ,f\left( t\right)
-m_{i}e_{i}\right\rangle \geq 0\text{ \ }  \label{3.10}
\end{equation}%
or, equivalently,%
\begin{equation*}
\left\Vert f\left( t\right) -\frac{M_{i}+m_{i}}{2}e_{i}\right\Vert \leq 
\frac{1}{2}\left( M_{i}-m_{i}\right) \text{ }
\end{equation*}%
\ for \ $i\in \left\{ 1,\dots ,n\right\} $ \ and a.e. $t\in \left[ a,b\right]
.$ Then we have the reverse of the generalised triangle inequality%
\begin{equation*}
\left[ \sum_{i=1}^{n}\frac{4m_{i}M_{i}}{\left( m_{i}+M_{i}\right) ^{2}}%
\right] ^{\frac{1}{2}}\int_{a}^{b}\left\Vert f\left( t\right) \right\Vert
dt\leq \left\Vert \int_{a}^{b}f\left( t\right) dt\right\Vert ,
\end{equation*}%
with equality if and only if%
\begin{equation*}
\int_{a}^{b}f\left( t\right) dt=\int_{a}^{b}\left\Vert f\left( t\right)
\right\Vert dt\left( \sum_{i=1}^{n}\frac{2\sqrt{m_{i}M_{i}}}{m_{i}+M_{i}}%
e_{i}\right) .
\end{equation*}
\end{corollary}

\begin{proof}
From the proof of Corollary \ref{c2.3}, we know (\ref{3.10}) implies that%
\begin{equation*}
\frac{2\sqrt{m_{i}M_{i}}}{m_{i}+M_{i}}\left\Vert f\left( t\right)
\right\Vert \leq \func{Re}\left\langle f\left( t\right) ,e_{i}\right\rangle
,\ \ \ i\in \left\{ 1,\dots ,n\right\} \text{ \ and a.e. }t\in \left[ a,b%
\right] .
\end{equation*}%
Now, applying Theorem \ref{t3.1} for $k_{i}:=\frac{2\sqrt{m_{i}M_{i}}}{%
m_{i}+M_{i}},$ $i\in \left\{ 1,\dots ,n\right\} ,$ we deduce the desired
result.
\end{proof}

\section{Applications for Complex-Valued Functions}

Let $e=\alpha +i\beta $ $\left( \alpha ,\beta \in \mathbb{R}\right) $ be a
complex number with the property that $\left\vert e\right\vert =1,$ i.e., $%
\alpha ^{2}+\beta ^{2}=1.$

The following proposition holds.

\begin{proposition}
\label{p4.1}If $f:\left[ a,b\right] \rightarrow \mathbb{C}$ is a Lebesgue
integrable function with the property that there exists a constant $K\geq 1$
such that%
\begin{equation}
\left\vert f\left( t\right) \right\vert \leq K\left[ \alpha \func{Re}f\left(
t\right) +\beta \func{Im}f\left( t\right) \right]  \label{4.1}
\end{equation}%
for a.e. $t\in \left[ a,b\right] ,$ where $\alpha ,\beta \in \mathbb{R}$, $%
\alpha ^{2}+\beta ^{2}=1$ are given, then we have the following reverse of
the continuous triangle inequality:%
\begin{equation}
\int_{a}^{b}\left\vert f\left( t\right) \right\vert dt\leq K\left\vert
\int_{a}^{b}f\left( t\right) dt\right\vert .  \label{4.2}
\end{equation}%
The case of equality holds in (\ref{2.2}) if and only if%
\begin{equation*}
\int_{a}^{b}f\left( t\right) dt=\frac{1}{K}\left( \alpha +i\beta \right)
\int_{a}^{b}\left\vert f\left( t\right) \right\vert dt.
\end{equation*}
\end{proposition}

The proof is obvious by Theorem \ref{t2.1}, and we omit the details.

\begin{remark}
\label{r4.2}If in the above Proposition \ref{p4.1} we choose $\alpha =1,$ $%
\beta =0,$ then the condition (\ref{4.1}) for $\func{Re}f\left( t\right) >0$
is equivalent to%
\begin{equation*}
\left[ \func{Re}f\left( t\right) \right] ^{2}+\left[ \func{Im}f\left(
t\right) \right] ^{2}\leq K^{2}\left[ \func{Re}f\left( t\right) \right] ^{2}
\end{equation*}%
or with the inequality:%
\begin{equation*}
\frac{\left\vert \func{Im}f\left( t\right) \right\vert }{\func{Re}f\left(
t\right) }\leq \sqrt{K^{2}-1}.
\end{equation*}%
Now, if we assume that%
\begin{equation}
\left\vert \arg f\left( t\right) \right\vert \leq \theta ,\ \ \ \theta \in
\left( 0,\frac{\pi }{2}\right) ,  \label{4.3}
\end{equation}%
then, for $\func{Re}f\left( t\right) >0,$%
\begin{equation*}
\left\vert \tan \left[ \arg f\left( t\right) \right] \right\vert =\frac{%
\left\vert \func{Im}f\left( t\right) \right\vert }{\func{Re}f\left( t\right) 
}\leq \tan \theta ,
\end{equation*}%
and if we choose $K=\frac{1}{\cos \theta }>1,$ then%
\begin{equation*}
\sqrt{K^{2}-1}=\tan \theta ,
\end{equation*}%
and by Proposition \ref{p4.1}, we deduce%
\begin{equation}
\cos \theta \int_{a}^{b}\left\vert f\left( t\right) \right\vert dt\leq
\left\vert \int_{a}^{b}f\left( t\right) dt\right\vert ,  \label{4.4}
\end{equation}%
which is exactly the Karamata inequality (\ref{1.2}) from the Introduction.
\end{remark}

Obviously, the result from Proposition \ref{p4.1} is more comprehensive
since for other values of $\left( \alpha ,\beta \right) \in \mathbb{R}^{2}$
with $\alpha ^{2}+\beta ^{2}=1$ we can get different sufficient conditions
for the function $f$ such that the inequality (\ref{4.2}) holds true.

A different sufficient condition in terms of complex disks is incorporated
in the following proposition.

\begin{proposition}
\label{p4.3}Let $e=\alpha +i\beta $ with $\alpha ^{2}+\beta ^{2}=1,$ $r\in
\left( 0,1\right) $ and $f:\left[ a,b\right] \rightarrow \mathbb{C}$ a
Lebesgue integrable function such that%
\begin{equation}
f\left( t\right) \in \bar{D}\left( e,r\right) :=\left\{ z\in \mathbb{C}|%
\text{ }\left\vert z-e\right\vert \leq r\right\} \text{ \ \ for a.e. }t\in %
\left[ a,b\right] .  \label{4.5}
\end{equation}%
Then we have the inequality%
\begin{equation}
\sqrt{1-r^{2}}\int_{a}^{b}\left\vert f\left( t\right) \right\vert dt\leq
\left\vert \int_{a}^{b}f\left( t\right) dt\right\vert .  \label{4.6}
\end{equation}%
The case of equality holds in (\ref{4.6}) if and only if%
\begin{equation*}
\int_{a}^{b}f\left( t\right) dt=\sqrt{1-r^{2}}\left( \alpha +i\beta \right)
\int_{a}^{b}\left\vert f\left( t\right) \right\vert dt.
\end{equation*}
\end{proposition}

The proof follows by Corollary \ref{c2.2} and we omit the details.

Finally, we may state the following proposition as well.

\begin{proposition}
\label{p4.4}Let $e=\alpha +i\beta $ with $\alpha ^{2}+\beta ^{2}=1$ and $%
M\geq m>0.$ If $f:\left[ a,b\right] \rightarrow \mathbb{C}$ is such that%
\begin{equation}
\func{Re}\left[ \left( Me-f\left( t\right) \right) \left( \overline{f\left(
t\right) }-m\overline{e}\right) \right] \geq 0\text{ \ \ for a.e. }t\in %
\left[ a,b\right] ,  \label{4.7}
\end{equation}%
or, equivalently,%
\begin{equation}
\left\vert f\left( t\right) -\frac{M+m}{2}e\right\vert \leq \frac{1}{2}%
\left( M-m\right) \text{ \ \ for a.e. }t\in \left[ a,b\right] ,  \label{4.8}
\end{equation}%
then we have the inequality%
\begin{equation}
\frac{2\sqrt{mM}}{M+m}\int_{a}^{b}\left\vert f\left( t\right) \right\vert
dt\leq \left\vert \int_{a}^{b}f\left( t\right) dt\right\vert ,  \label{4.9}
\end{equation}%
or, equivalently,%
\begin{equation}
\left( 0\leq \right) \int_{a}^{b}\left\vert f\left( t\right) \right\vert
dt-\left\vert \int_{a}^{b}f\left( t\right) dt\right\vert \leq \frac{\left( 
\sqrt{M}-\sqrt{m}\right) ^{2}}{M+m}\left\vert \int_{a}^{b}f\left( t\right)
dt\right\vert .  \label{4.10}
\end{equation}%
The equality holds in (\ref{4.9}) (or in the second part of (\ref{4.10})) if
and only if%
\begin{equation*}
\int_{a}^{b}f\left( t\right) dt=\frac{2\sqrt{mM}}{M+m}\left( \alpha +i\beta
\right) \int_{a}^{b}\left\vert f\left( t\right) \right\vert dt.
\end{equation*}
\end{proposition}

The proof follows by Corollary \ref{c2.3} and we omit the details.

\begin{remark}
\label{r4.5}Since%
\begin{align*}
Me-f\left( t\right) & =M\alpha -\func{Re}f\left( t\right) +i\left[ M\beta -%
\func{Im}f\left( t\right) \right] , \\
\overline{f\left( t\right) }-m\overline{e}& =\func{Re}f\left( t\right)
-m\alpha -i\left[ \func{Im}f\left( t\right) -m\beta \right] 
\end{align*}%
hence%
\begin{multline}
\func{Re}\left[ \left( Me-f\left( t\right) \right) \left( \overline{f\left(
t\right) }-m\overline{e}\right) \right]   \label{4.11} \\
=\left[ M\alpha -\func{Re}f\left( t\right) \right] \left[ \func{Re}f\left(
t\right) -m\alpha \right] +\left[ M\beta -\func{Im}f\left( t\right) \right] %
\left[ \func{Im}f\left( t\right) -m\beta \right] .
\end{multline}%
It is obvious that, if%
\begin{equation}
m\alpha \leq \func{Re}f\left( t\right) \leq M\alpha \text{ \ \ for a.e. }%
t\in \left[ a,b\right] ,  \label{4.12}
\end{equation}%
and 
\begin{equation}
m\beta \leq \func{Im}f\left( t\right) \leq M\beta \text{ \ \ for a.e. }t\in %
\left[ a,b\right] ,  \label{4.13}
\end{equation}%
then, by (\ref{4.11}),%
\begin{equation*}
\func{Re}\left[ \left( Me-f\left( t\right) \right) \left( \overline{f\left(
t\right) }-m\overline{e}\right) \right] \geq 0\text{ \ \ for a.e. }t\in %
\left[ a,b\right] ,
\end{equation*}%
and then either (\ref{4.9}) or (\ref{4.12}) hold true.
\end{remark}

We observe that the conditions (\ref{4.12}) and (\ref{4.13}) are very easy
to verify in practice and may be useful in various applications where
reverses of the continuous triangle inequality are required.

\begin{remark}
Similar results may be stated for functions $f:\left[ a,b\right] \rightarrow 
\mathbb{R}^{n}$ or $f:\left[ a,b\right] \rightarrow H,$ with $H$ particular
instances of Hilbert spaces of significance in applications, but we leave
them to the interested reader.
\end{remark}

\end{document}